\documentclass[11pt,bezier]{article}
\usepackage{amsmath,amssymb,amsfonts,euscript}
\usepackage[usenames]{color}
\newpage

\newcommand{\be}{\begin{equation}}
\newcommand{\ee}{\end{equation}}
\newtheorem{theorem}{Theorem}[section]
\newtheorem{lemma}{Lemma}[section]

\newtheorem{definition}{Definition}[section]
\newtheorem{remark}{Remark}[section]
\newtheorem{proposition}{Proposition}[section]
\renewcommand{\theequation}{\arabic{section}.\arabic{equation}}
\title{\bf\Large Certain Periodically Correlated Multi-component Locally Stationary Processes}
\vspace{-1cm}
\author
{ N. Modaresi\,\,\,\, and \,\,\,  S. Rezakhah\thanks{\scriptsize
Faculty of Mathematics and Computer Science, Amirkabir University of
Technology, 424 Hafez Avenue, Tehran 15914, Iran. E-mail:
namomath@aut.ac.ir(N. Modaresi), rezakhah@aut.ac.ir(S. Rezakhah).}
}
\date{}
\begin{document}

\maketitle

\begin{abstract}
By introducing $X^{ls}(t)$ as a random mixture of two stationary processes where the time dependent random weights have  exponentially convex covariance, we show that this process has a  multi-component locally stationary covariance function in Silverman's sense.
We also define $X^p(t)$ as a certain continuous time periodically correlated (PC) process where its covariance function is generated by the covariance function of a discrete time through defining some simple random measure on real line.  We also impose a bi-periodic correlation for this PC process with $X^{ls}(t)$.
The existence of such  random measure is proved.
Then by defining  $X(t)=X^{ls}(t)+X^p(t)$ as  a certain periodically correlated multi-component locally stationary process,
the covariance structure and time varying spectral representation of such  processes are characterized.\\

\end{abstract}

\section{Introduction}
Recently a large amount of work has been devoted to time series analysis, with the focus placed on locally stationary (LS) and periodically correlated (PC) processes which give plausible description of real world.
Silverman \cite {w1} presented a definition of LS process to model systems which behave as a function of time. He presented a relation between the covariance and the spectral density of LS processes which constitutes a natural generalization of the Wiener-Khintchine relation.
A general class of LS processes where their spectral structure varies smoothly over time is introduced by Dahlhaus \cite{d0}.
A computationally efficient state space method for estimating, predicting and making statistical inferences about these non-stationary models are proposed in  \cite{p1}.
Exponentially convex stochastic processes, their covariance functions and spectral representation have been studied by Loeve \cite {l1} among others \cite{e1}.

Gladyshev \cite {g1} introduced PC sequences and showed that all are harmonizable in the sense of Loeve \cite {l1}. He provided the structure of the covariance function and an interesting spectral representation.
At the most basic level, the connection of the PC structure to a group of shift and unitary operators was mentioned as a motivation for the development of spectral representations that are much like those of Gladyshev.
One of the classical results of unitary operator theory is the spectral theorem, the notion of spectral measures and the time varying spectral representation of a PC process \cite{h1}.
The contributions to the analysis of PC processes in the literature are \cite{d2}, \cite{h1} and \cite{s1}, where various properties of PC processes are presented.
Harmonic series representation, coherent and component method for mean and covariance function estimation by using linear filtration
 of PC processes  are studied in \cite{j1} and \cite{j2}.
A class of bilinear processes with periodic time-varying coefficients of periodic ARMA and periodic GARCH models has been applied in \cite{b1}.

We consider a certain non-stationary process $\{X(t), t\in {\mathbb R}^+\}$, say  periodically correlated multi-component locally stationary (PC-LS) process.
For introducing the structure of this process we consider a partition of the positive real line.
{\bf This paper can be considered as an extension of previous works in three stages.
First by presenting spectral representation of certain continuous PC process. This is established by introducing some random measure on subsets of some proper partition of real line where
the increments of the process on such partition, provides a discrete time PC process.
 We specify the covariance structure of such random measure through the covariance structure of the discrete time PC process, and provide the spectral representation of such continuous time PC process.  The intuition of this process comes from the aggregated traffic on various networks.
Second we present a locally stationary process in the Silverman sense as a time varying random mixture of stationary processes where the coefficients have exponentially convex covariances and are independent of the stationary ones.
This consideration provides smooth transition in the correlation structure of the process which is in effect of two or more stationary resources at the nodes of the mentioned partitions. This makes a convenient way for analyzing traffic systems like the mobile cellular networks.
Third we introduce the process as sum of such PC and LS processes which have bi-periodic correlations and have potential to model periodic and local variation 
in network traffics.
The idea of introducing such process comes from the behavior of traffic flow in different networks, like telecommunication, internet or transportation networks.} This process is determined by aggregation of the effects of two main components. The first component often has periodic behavior of the flow in time, $X^p(t)$, which shows over all usage of the network. We assume that $X^p(t)$ in turn can be approximated via a periodic sequence $X_j^p$, which represents the increment  of the flow on some proper partition of time, say $B_j, j\in {\Bbb N}$. Type and amount of the usage of the network by different users at a time have some interaction which effects the flow and provide the other component of the process and can be approximated as a  multi-component locally stationary process $X^{ls}(t)$.
We assume that  $X^{ls}(t)$ is  a random mixture, with some exponentially convex random weight process, of two stationary processes at time points inside each $B_j$.
We consider some fixed $T\in {\mathbb N}$ and a partition of the positive real line as $0=s_0<s_1<\ldots$,
$\, B_j=(s_{j-1}, s_j]$ for $j\in\mathbb{N}$ and  $|B_j|=|B_{j+kT}|$, $k\in\mathbb{N}$, then $S=\sum_{i=1}^T|B_i|$. Let $X(t)=X^{ls}(t)+X^p(t)$, $t\in{\mathbb{R}^+}$ represents a stochastic process, where $X^{ls}(t)=\sum_{j=1}^{\infty}X^{ls}_j(t)I_{B_j}(t)$ and $X^{ls}_j(t)$ is a multi-component LS process which is an exponentially convex mixture of two stationary processes.
Also let $\{X^p_j\}$, $j\in\mathbb{N}$  be the discrete time PC process, $M_j$, $j\in\mathbb{N}$ the sequence
of random measure on $B_j$ where $M_j(B_j):=X^p_j$ and $X^p_j(t):=M_j(s_{j-1}, t]$ for $t\in B_j, \; j\in{\mathbb N}$ with some special correlation. As the amount of the process on a set  $A$, say traffic, we consider simple random measure
$M(A)=\sum_{j=1}^m M_{i_j}(A\cap B_{i_j})$ where $\;A=\cup_{j=1}^m (A\cap B_{i_j})$. We prove the existence of such measures by introducing some measurable map and some $T$-variate random measure.
Then we study $X^p(t)=\sum_{j=1}^{\infty}\Delta^p_j(t)I_{B_j}(t)$, where $\Delta^p_j$ is a linear combination of $X^p_{j-1}$ and $X^p_j(t)$, which describes elimination of the effect of traffic on $B_{j-1}$ as a decreasing fraction of $X^P_{j-1}$ and aggregate of the process on $B_j$ as $X^p_j(t)$, and study its spectral representation.

In section 2, we study the general framework and preliminaries of the harmonizable representation of stationary an PC processes and also, we present exponentially convex processes. The main result as the covariance structure and spectral representation of such continuous time multi-component PC-LS processes are given in section 3.

\renewcommand{\theequation}{\arabic{section}.\arabic{equation}}
\section{Preliminaries}
\setcounter{equation}{0}
We review the harmonizable representations of PC processes based on unitary operators and also the definition of exponentially convex processes in this section.

\subsection{Spectral representation of PC processes}
We review harmonizable representation of PC processes by unitary operators.  For a comprehensive  review of these one can refer to  Hurd and Miamee \cite {h1}.
Every wide sense stationary process $X(t)$ has harmonizable representation
\vspace{-1mm}
\be X(t)=\int_{-\infty}^{\infty}e^{i\lambda t}dZ(\lambda)\label{sbn0-0}\ee
where $Z(\lambda)$ has orthogonal increments.

A unitary operator on a Hilbert space $\cal{H}$ is a linear operator from $\cal{H}$ to $\cal{H}$ which preserve inner product as $\big< Ux,Uy\big>=\big<x,y\big>=\mbox{cov}(x, y)$ for every $x,y \in \cal{H}$.
\begin{theorem}
For any unitary operator $U$ on a Hilbert space $\cal{H}$, there exists a unique spectral measure $Q$ on the Borel subsets of $[0,2\pi)$ such that $U=\int_{0}^{2\pi} e^{i\lambda} Q(d\lambda)$, and $U^t=\int_{0}^{2\pi} e^{i\lambda t} Q(d\lambda)$ for any integer $t$.
\end{theorem}
Existence of unitary operator for a PC sequence, characterize its spectral representation.
\begin{proposition}
A second order stochastic sequence $X_j$ is PC with period $T$ if and only if for every $j\in\mathbb{Z}$, there exists a unitary operator $U=V^T$ and a periodic sequence (process) $P_j$ with period $T$ taking values in ${\cal{H}}_X=\overline{sp}\{X_j, j\in\mathbb{Z}\}$ for which
$X_j=V^j P_j$
where $V=\int_{0}^{2\pi} e^{i\lambda/T}Q(d\lambda)$ and $Q$ is the spectral measure defined in Theorem $2.1$, so
$X_j=\int_{0}^{2\pi}e^{i\lambda j/T}Q(d\lambda)P_j.$
\end{proposition}

\subsection{Exponentially convex process}
We give a brief description of exponentially convex process and its covariance function, for more details see \cite {e1}, \cite {g3}.

\begin{definition}
The covariance function of a second order zero mean process $\{Z(t), t\in\mathbb{R}\}$ with finite variance,
$\big<Z(t_i),\overline{Z(t_j)}\big>=\psi(t_i+t_j)$, is called exponentially convex if and only if
$\sum_{i=1}^{n}\sum_{j=1}^{n}a_i\overline{a_j}\psi(t_i+t_j)\geqslant 0$,
for all finite sets of complex coefficients $a_1, \ldots, a_n$ and points $t_1, \ldots, t_n\in\mathbb{R}$.
\end{definition}
 A stochastic process with such covariance function is called exponentially convex. The result of Berg et al. \cite {b0} implies that a continuous function is exponentially convex if and only if it is the Laplace transform of a non-negative finite measure.

\renewcommand{\theequation}{\arabic{section}.\arabic{equation}}
\section{Main results: A continuous time PC-LS model}
\setcounter{equation}{0}
We introduce a new class of certain non-stationary process, say multi-component periodically correlated locally stationary (PC-LS) process as
\be X(t)=X^{ls}(t)+X^p(t),\hspace{7mm}t>0\label{sbn0-1}\ee
where $X^p(\cdot)$ is a continuous time PC process and $X^{ls}(\cdot)$ is a multi-component LS process which have bi-periodic correlation.  We introduce their special structures  in subsections 3.1 and 3.2 respectively. Their covariance and cross covariance structures are studied in subsection 3.3. Finally the spectral representation of the process is characterized in 3.4.

\subsection{Continuous PC process}
Let $\{X^p_j, j\in\Bbb{N}\}$ be a positive second order discrete time PC process with period $T$,
$M_j\,$ a random measure on Borel field of subsets of $B_j$, where $M_j(B_j):=X^p_j$, and for $A, B\subset B_j$,
$\; D\subset B_k$,  we define $\;E[M_j(A)]=\frac{|A|}{|B_j|}E[X^p_j]$  and covariance functions \vspace{-1mm}
\be {\bf\gamma_{j,j}(A,B)\!=\!\frac{\lambda |A||B|\!+\!(1\!-\!\lambda )a_j|A\cap B|}{a_j^2}\gamma^p_{jj}}, \hspace{.1in} \gamma_{j,k}(A,D)\!:=\!\frac{|A||D|}{a_ja_k}\gamma^p_{jk} \ee \vspace{-3mm}

\noindent
where $\gamma_{j,j}(A,B)=\big<M_j(A), M_j(B) \big>$,  $\gamma_{j,k}(A,D)=\big<M_j(A), M_k(D) \big>$, $\gamma^p_{jk}=\big< X^p_j, X^p_k\big>$,  $\gamma^p_{jj}=\mbox{Var} (X^p_j)$,  $0\leq \lambda \leq 1 \;$ and $a_j=|B_j|$.\\
If $B=B_k$ then for $j\neq k$,
$\big<M_j(A),M_k(B_k)\big>=\frac{|A|}{a_j}\gamma^p_{jk}.$
This inner product is well defined.
Let $X^p_j(t):=M_j(s_{j-1}, t]$, $t\in B_j$, $j\in{\mathbb N}$. So for $t,u\in B_j, t\leq u$, and  $v\in B_k$
\vspace{-.1mm}
\be\big<X^p_j(t),X^p_j(u)\big>=\frac{\lambda \, a^t_ja^u_j+(1-\lambda )a_j\, a^t_j}
{a_j^2}\gamma^p_{jj},\;\;\;\;\;\big<X^p_j(t),X^p_k(v)\big>
=\frac{a^t_ja^v_k}{a_ja_k}\gamma^p_{jk}\label{sbn3}\ee \vspace{-3mm}

\noindent
where $a^t_j=t-s_{j-1}$. Thus by defining
\vspace{-1mm}

\be N_j(t-s_{j-1}):=M_j(s_{j-1},t]=X^p_j(t),\label{sbn2}\ee \vspace{-4mm}

\noindent
we find that  $N_j(y)$ is a discrete time PC process with period $T$ with respect to $j$ for fixed $0< y\leqslant a_j$ and $X_j^p(t)$ is a bi-periodic process with period $T$ in $j$ and period $S=\sum_{i=1}^T |B_i|$ in $t$, that is  $<X_j^p(t),X_k^p(u)>=<X_{j+T}^p(t+S), X_{k+T}^p(u+S)>$. Also for $t>0$
\vspace{-2mm}
\be{\bf X^p(t)=\sum_{j=1}^{\infty}\Delta^p_j(t)I_{B_j}(t);\hspace{5mm}\Delta_j^p(t)
=\frac{a_j-a_j^t}{a_j}X^p_{j-1}+X_j^p(t),\label{sbn0-3}}\ee \vspace{-4mm}

\noindent
is a continuous time PC process with respect to $t$ with period  $S$.
{\bf Using a similar method such as the one described by Soltani and Parvardeh \cite{s2}, we prove the existence of such random measures by the following.

Let $(D, {\cal D})$ be a measurable space and ${\cal L}^2(\Omega,{\cal F},P)$ the Hilbert space of real random variables on the probability space $(\Omega, {\cal F},P)$ with finite second moment. A mapping $\Phi: {\cal D} \rightarrow {\cal L}^2(\Omega, {\cal F},P)$ is a second order random measure on ${\cal D}$ if $\Phi(\emptyset  )=0$ and $E\big|\Phi \big(\cup_{i=1}^{\infty} A_i\big)-\sum_{i=1}^n \Phi (A_i)\big|^2\rightarrow 0$ as $n\rightarrow \infty$, for disjoint sets
$A_1, A_2, \ldots  \in {\cal D}$. We introduce simple random measure $M(A)$ for $A\in {\cal D}$ by
$M(A)=M_1(A\cap B_1)+M_2(A\cap B_2)+\ldots +M_T(A\cap B_T)$ where $B_j$ is the support of $M_j,\; j=1,\ldots, T$ and
$B_i\cap B_j= \emptyset $ for $i \neq j$. For the random measure $M$, $E|M(dx)|^2=\gamma_{j,j}(dx)$ for $x \in B_j$. Also for $x\in B_j$ and $y\in B_i$ we have $E\big(M(dx)M(dy)\big)=E\big(M_j(dx), M_i(dy)\big)=\gamma_{j,k}(dx, dy)$. To show the existence of such random measures, let $\Psi=(\Psi_1,\ldots \Psi_T)$ be a $T$-variate random measure on $(X,{\cal X})$, where $X=B_1\times \ldots \times B_T$ and ${\cal X}$ is the corresponding class of finite disjoint union of measurable rectangles.
In this case for $A\in {\cal D}$, $A=\cup_{j=1}^m (A\cap B_{i_j})$ and  $\{i_1,i_2,\ldots ,i_m\}$ are $m$ distinct integer of $\{1,2,\ldots, T\}$. Let $X^*=B_{i_1}\times \ldots \times B_{i_m}$ and ${\cal X}^*$ the class of finite disjoint rectangles $F^*=(F\cap B_{i_1})\times \ldots \times (F\cap B_{i_m})$ where $F\subset D^*=\cup_{j=1}^m B_{i_j}$ and $S:{\cal X}^*\rightarrow  D^*$ be a measurable map that $S(A^*)=\cup_{i=1}^m \{x_i:\underline{x}=(x_1,\ldots , x_m)\in A^*\}$ for $A^*\in {\cal X}^*$.
Then $S^{-1}(F)=\{\underline{x}=(x_1,\cdots , x_m): x_j\in F\cap B_{i_j}, j=1,\ldots, m\}$ for $F\subset { D}^*$. Also let
$M_{i_j}(F)=\Psi_{i_j}S^{-1}(F)$, and $\Psi^*=\big(\Psi_{i_1} ,\ldots, \Psi_{i_m}\big) $ be an $m$-variate random measure on $(X^*,{\cal X}^*)$  with covariance matrix {\Large ${\Bbb\mu}\;$}$=[\mu_{j,k}]$, then for the resulting $M$ and for $1\leq j,k\leq m$
$$
\mu_{j,k}\big(S^{-1}(F)\big)=E\big(\Psi_{i_j}(S^{-1} (F))\Psi_{i_k}(S^{-1}(F))\big)=\gamma_{i_j,i_k}(F\cap B_{i_j},F\cap B_{i_k})
$$
$$
\mu_{j,j}\big(S^{-1}(F)\big)=E\big(\Psi_{i_j}(S^{-1} (F))\Psi_{i_j}(S^{-1}(F))\big)=\gamma_{i_j,i_j}(F\cap B_{i_j}, F\cap B_{i_j})
$$
Interestingly
$$
\nu(A\times C)=\sum_{j,k} \gamma_{j,k}(A\cap B_j, C\cap B_k)
$$
defines a product measure on $(D^*\times D^*, {\cal Y})$, where ${\cal Y}$ is the class of finite disjoint union of corresponding measurable rectangles. Also for $E_{j,k}\in {\cal Y}$ where $E_{j,k}\subset B_{i_j}\times B_{i_k}$
$$
\nu (E_{j,k})=\mu_{j,k}\big\{\underline{x}: \underline{x}\in X^*, (x_j,x_k)\in E_{j,k}\big\}
$$
Thus for $E\in{\cal Y}$ and $E_{j,k}=E\cap (B_{i_j}\times B_{i_k})$
\begin{eqnarray*}
\nu (E)&\!\!=\!\!&\sum_{j,k} \nu(E_{j,k})=\sum_{l=-m+1}^{m-1} \sum_{j=\max \{1-l, 1\}}^{\min \{m-l,m\}} \nu (E_{j,l+j})\\&=&\sum_{j=1}^m \nu(E_{j,j})+\sum_{l=1}^{m-1}\bigg\{\sum_{j=1}^{m-l}\nu(E_{j,l+j})
+\sum_{j=m-l+1}^{\min \{2m-l,m \}}\nu (E_{j,l+j-m})\bigg\}.
\end{eqnarray*}
Now define $\mu_0,\ldots, \mu_{m-1}$ on $D^*$ through
\begin{eqnarray*}
\mu_0(A)& = &\sum_{j=1}^m \mu_{j,j}\big(S^{-1}(A)\big), \\
\mu_l(A)& = &\sum_{j=m-l+1}^{\min \{2m-l,m\}}\mu_{j,j+l-m}\big(S^{-1}(A) \big)
+\sum_{j=1}^{m-l}\mu_{j,j+l}\big(S^{-1}(A)) \big)
\end{eqnarray*}
for $l=1,\ldots, m-1$; then  for $E\in {\cal Y}$ we have that
$\nu (E)=\sum _{l=0}^{m-1} \mu_l \{ x \in D^*; (x,y) \in E \;\; \mbox{for some }\; y \}.$
Therefore the product measure $\nu$ is specified by $m$ set functions $\mu_0, \ldots, \mu_{m-1}$ on $D^*$ that are determined by {\Large$ \mu$}$=[\mu_{j,k}]$, the covariance matrix of $m$-variate random vector $M=(M_1, \ldots, M_m)$.
  The covariance matrix       {\Large $\mu$} also can be specified from $\mu_0, \ldots, \mu_{m-1}$, namely for $G\in {\cal X}^*$
\begin{eqnarray*}
\mu_{j,k}(G)&\!\!=\!\!&\mu_{m+k-j}\big(\{x_k: { \underline{x}\in G}\}\cup \{x_j: { \underline{x}\in G}\} \big),\hspace{1cm} k-j<0, \\
\mu_{j,k}(G)&\!\!=\!\! & \mu_{k-j}\big(\{x_k: { \underline{x}\in G}\}\cup \{x_j: { \underline{x}\in G}\} \big),\hspace{15mm} k-j\geq 0,
\end{eqnarray*}
Therefore  $ \nu$ and {\Large $\mu$} uniquely specify each other.}

\subsection{Multi-component LS process}
{\bf To provide the smoothing transition and LS behavior of the process between inside each partition we introduce a random mixture of two or more stationary processes or resources at nodes of such partition with exponentially convex coefficients which provides an LS processes in the Silverman sense.}
Let $\big \{Y^s_j(t), t\in {B_{j}\cup B_{j+1}}\big \}$,
$Y^s_j(t)=\int_{-\infty}^{\infty}e^{it\lambda}\eta_j(d\lambda)$, $j\in\mathbb{N}$ be a sequence of independent stationary processes, $\eta_j(\cdot)\stackrel{d}{=}\eta_{j+T}(\cdot)$
and $\{Y_j^s(t)\}\stackrel{d}{=}\{Y_{j+kT}^s(t+kS)\}$ for $k\in\Bbb N$. We define
\be X^{ls}(t)=\sum_{j=1}^{\infty}X^{ls}_j(t)I_{B_j}(t),\hspace{5mm}X^{ls}_j(t)\!=\!U^{j-1}(t)Y^s_{j-1}(t)+U^j(t)Y^s_j(t)\label{sbn0-3-1}\ee
for $t\in B_j$, $j=2, \cdots, T$ and $X_1^{ls}(t)=U^1(t)Y_1^s(t)$ for $t\in B_1$, in which $\{U^j(t), j\in \mathbb{N} \}$ is a random weight process with exponentially convex covariance and are independent of the process  $Y^s_j(\cdot)$. We call $X_j^{ls}(\cdot)$, a multi-component LS process motivated from its covariance function, which is represented in Lemma 3.1.
%

\subsection{Covariance function of PC-LS process}
We present some lemmas in this section which we use to provide the covariance function of the process $X(t)$ in Theorem 3.1.

\begin{lemma}
Let $B_j$, $j\in\mathbb{N}$ be the introduced partition of positive real line. The covariance function of $X^{ls}(t), t\in {\mathbb R}^+$ satisfies
\vspace{-3mm}
\be\gamma^{ls}(t,u)\equiv\big<X^{ls}(t), X^{ls}(u)\big>
=\sum_{m=1}^{\infty}\sum_{n=m-1}^{m+1}\gamma^{ls}_{mn}(t,u)I_{B_m}(t)I_{B_n}(u)\label{sbn3-2}\ee
where $\;\gamma^{ls}_{mm}(t,u)\equiv \big<X_m^{ls}(t),X_m^{ls}(u) \big>=\psi_{m}(t+u)\gamma_{m}(t-u)+ \psi_{m-1}(t+u)\gamma_{m-1}(t-u),\,$
$\, \gamma^{ls}_{mn}(t,u)\equiv \big<X_m^{ls}(t),X_n^{ls}(u) \big>=\psi_{k}(t+u)\gamma_{k}(t-u)$, $k=\min\{m,n\}$ and $|n-m|=1$ in which
$\psi_m(t+u)=\big<U^m(t),U^m(u)\big>$, $\gamma_m(t-u)=\big<Y^s_m(t),Y^s_m(u)\big >$ for $m\in \mathbb{N}$.
\end{lemma}
{\bf Proof:}
According to $(\ref{sbn0-3-1})$ and the fact that $U^j(\cdot)$ and $Y^s_j(\cdot)$ are independent processes,
so $\gamma^{ls}_{mm}(t,u)$ and $\gamma^{ls}_{mn}(t,u)$ for $|n-m|=1$ are as presented by the lemma.
Thus $X^{ls}_j(\cdot)$ is a multi-component LS process in the Silverman sense \cite{w1}. Also by (3.6)
$$\gamma^{ls}(t,u)=\sum_{m=1}^{\infty}\big<X^{ls}_m(t), X^{ls}(u)\big>I_{B_m}(t)
=\sum_{m=1}^{\infty}\sum_{n=m-1}^{m+1}\gamma^{ls}_{mn}(t,u)I_{B_m}(t)I_{B_n}(u).$$

\begin{lemma}
 Let $B_j$, $j\in\mathbb{N}$ be a partition of positive real line and $\gamma^p_{mn}=\big<X^p_m,X^p_n\big>$. The covariance function of $X^p(t)$ for $t\in B_m$ and $u\in B_n$, $t\leq u$ is
\be\hspace{-1in}\gamma^p(t,u)\equiv\big<X^p(t), X^p(u)\big>= \hspace{3in}\label{sbn3-3}\ee
$$\left\{\begin
{array}{cccc}
\vspace{-1mm}&& \hspace{-.8cm}
{\bf A_{t,m}\big(A_{u,m}\gamma_{m-1,m-1}\!+\!\frac{a_m^u}{a_m}\gamma_{m-1,m}\big)\!+\!A_{u,m}\frac{a_m^t}{a_m}\gamma_{m-1,m}\!+
\!B_m(t,u)\gamma^p_{mm} }\hspace{5mm} n=m,\vspace{2mm} \\&& \hspace{-1cm}
{\bf A_{t,m}\big( A_{u,n} \gamma_{m-1.n-1}^p\!+ \!\frac{a_n^u}{a_n}\gamma_{m-1,n}^p\big)\!\!+\!A_{u,n}\frac{a_m^t}{a_m}\gamma_{m,n-1}^p\!+\!\frac{a_m^ta_n^u}{a_ma_n}\gamma_{m,n}^p } \hspace{1cm} \hspace{-3mm} |m-n|\geq 1\end {array}\right.\\$$
where $A_{t,m}=(1-a_m^t/a_m)$, $\;A_{u,n}=(1-a_n^u/a_n)$ and $B_m(t,u)=\!\frac{\lambda a^t_ma^u_m+(1-\lambda )a_m a^t_m}{a_m^2}$.
\end{lemma}
{\bf Proof:}
By (\ref{sbn0-3}) the covariance function of $X^p(\cdot)$  for $t\in B_m, u\in B_n$ is
$\gamma^p(t,u)=\big< X_m^p(t), X_n^p(u) \big>$ and by (\ref{sbn3})
we have the result.\\

\noindent
Let $a_j=|B_j|$, $a=E[U^j(t)]$ for $t\in {\mathbb R}^+$ and $\theta_{r}(j,u)=\big<X_j^p, Y^s_r(u)\big>$. We also define the cross covariance function of $X_j^p(t)=M_j(s_{j-1},t]$ and $Y^s_r(u)$  for $t\in B_j, \; u\in B_r$  as
\be\big<X_j^p(t), Y_r^s(u)\big>=\frac{a^t_j}{a_j}\theta_{r}(j,u),\;\;\; \big<X^{p}_j(t), X^{ls}_r(u)\big>=\frac{aa^t_j}{a_j}\big[\theta_{r-1}(j,u)+\theta_{r}(j,u)\big].   \ee

\begin{remark}
The cross covariance function of  $X^{ls}(t)$ and $X^p(t)$ can be written as
\be {\bf\gamma^{p,ls}(t,u)\equiv\big<X^{p}(t), X^{ls}(u)\big>
=a\sum_{m=1}^{\infty}\sum_{n=1}^{\infty}\big[ \frac{a^t_m}{a_m}D_{m,n}(u) +A_{t,m}D_{m-1,n}(u)\big]I_{B_m}(t)I_{B_n}(u) \label{sbn3-4}}\ee
where $D_{m,n}(u)=\theta_{n-1}(m,u)+\theta_{n}(m,u)$, $\; A_{t,m}=1-a_m^t/a_m\;$ and
$\;\theta_{r}(m,u)$ is defined by (3.9).
\end{remark}

\begin{theorem}
The covariance function of the multi-component {\em PC-LS} process $X(t)=X^{ls}(t)+X^p(t)$, where $X^{ls}(t)$ and $X^p(t)$ are dependent and defined by
$(\ref{sbn0-3-1})$ and $(\ref{sbn0-3})$ respectively, is
$\gamma(t,u)=\mbox{\em{cov}}\big(X(t), X(u)\big)=\gamma^{ls}(t,u)+\gamma^{p}(t,u)+\gamma^{p,ls}(t,u)$ where $\gamma^{ls}(t,u)$, $\gamma^{p}(t,u)$ and $\gamma^{p,ls}(t,u)$ are as in Lemma 3.1, Lemma 3.2 and Remark 3.1 respectively.
\end{theorem}

\subsection{Spectral representation}
In this section we obtain spectral representations of PC process $\{X^p(t), t\in {\mathbb R}^+\}$, and   locally stationary process  $\{X^{ls}(t), t \in {\mathbb R}^+\}$. Then by imposing bi-periodic property for the cross covariance function of $X^{ls}(t)$ and $X^p(t)$ their cross spectrum is characterized. Using these results the spectral representation of the multi-component PC-LS process $X(\cdot)$ and its spectral measure is provided by Theorem 3.2.

Let $\{X_j^p\}$, $j\in {\mathbb N}\}$ be a PC process with period $T$. By using {\bf the} Gladyshev result for PC processes, it has been shown in [9] that for $j,k \in {\mathbb N}$
\vspace{-2mm}
\be X_j^p=\int_0^{2\pi} e^{i \lambda j} {\bf d\vartheta(\lambda )},\label{vart-spec}\ee
where ${\bf d\vartheta(\lambda )=\sum_{k=0}^{T-1} Q(Td\lambda-2k\pi)\widetilde{P}_kI_{\Delta_k}(\lambda)}$, $\Delta_k=[2k\pi/T,2(k+1)\pi/T]$, $\widetilde{P}_k$ are Fourier coefficients of the periodic sequence $P_j$, $Q(\cdot)$ and $P_j$ are defined by Proposition 2.1, that $P_j=\sum_{k=0}^{T-1} \widetilde{P}_ke^{i2\pi k j/T}$.
The support of the spectral density
$\tilde{\theta}(d\lambda,d\omega)=\big<\vartheta(d\lambda),\vartheta(d\omega)\big>$
is the intersection of the set $\Psi =\{(\lambda,\omega): \lambda=\omega-2\pi k/T, k\in [-(T-1),T-1] \}$ and the square $[0,2\pi)\times [0,2\pi)$.

\begin{remark}
Let $\{X^p_j\}$ be a sequence of {\em PC} process with period $T$, by proposition (2.1)  $X^p_j=\int_0^{2\pi}e^{i\lambda j/T}\xi(d\lambda,j)$ where $\xi(d\lambda,j)=Q(d\lambda)P_j$,  $Q(d\lambda)$ is an orthogonally scattered random measure and $P_j$  a periodic sequence with period $T$. Let $a_j=|B_j|$, $a^t_j=t-s_{j-1}$  for $t\in B_j$. As by (3.4),
$X_j^p(t)\equiv N_j(a^t_j)$ and $N_j(y)$ is a discrete PC process with period $T$ for fixed $y$, so by Proposition $2.1$
$\; X^p_j(t)=\int_0^{2\pi}e^{i\lambda j/T}\zeta_j(d\lambda, t)$ for $t\in B_j$, where
$\zeta_j(d\lambda, t)=Q(d\lambda)\widehat{P}_{j,a^t_j}$
and $\widehat{P}_{j, y}$ is a periodic function in $j$ with period $T$ for fixed $y$. Using the relation (3.4) we have $X_j^p\equiv X_j^p(s_j)$, so $\widehat{P}_{j, a_j}\equiv P_{j}$. By a similar method as in  (3.11)
$\; X^p_j(t)=\int_0^{2\pi}e^{i\lambda j}\Upsilon (d\lambda, a^t_j)$ for
$t\in B_j$, where $\Upsilon (d\lambda, a^t_j)$ has the same structure as  $\vartheta(d \lambda )$ in (3.12),
by replacing $\widetilde{P}_k$ with $\widetilde{P}_{k,a_k^t}$ where $\widehat{P}_{j,a^t_j}=\sum_{k=0}^{T-1} \widetilde{P}_{k,a_k^t}e^{i2\pi j k/T}$. So $(\ref{sbn3})$ implies that
$$\Theta_{j,k}(d\lambda,d\omega,a_j^t,a_k^u)=\big<\Upsilon(d\lambda,a_j^t),\Upsilon(d\omega,a_k^u)\big>=\left\{\begin
{array}{cc}
\hspace{-3mm}{\bf \frac{\lambda a^t_ja^u_j+(1-\lambda )a_j a^t_j}{a_j^2} \tilde{\theta}(d\lambda,d\omega)} & j=k\\
\hspace{-3mm}\frac{a^t_ja^u_k}{a_ja_k}\, \tilde{\theta}(d\lambda,d\omega) & j\neq k\\
\end {array}\right.\hspace{1cm}$$
which is bi-periodic with period $T$ in $j$ and $k$, and period $S$ in $t$ and $u$ ($t\leq u$)
where as in (3.11),  $\tilde{\theta}(d\lambda,d\omega)=\big<\vartheta(d\lambda),\vartheta(d\omega)\big>$ have $\Psi$ as its support.
\end{remark}
{\bf We represent the spectral representation for $X^{ls}_j(\cdot)$ by the following remark.}

\begin{remark}
Spectral representation of the process $X^{ls}_j(\cdot)$, $j\in\mathbb{N}$ defined by $(\ref{sbn0-3-1})$ is
\be X^{ls}_j(t)=\int_{-\infty}^{\infty}e^{i\lambda t}\Phi_j(d\lambda, t),\hspace{7mm}t\in B_j\label{sbn5}\ee
where
$\Phi_j(d\lambda,t)=U^{j-1}(t)\eta_{j-1}(d\lambda)+ U^j(t)\eta_{j}(d\lambda)   $,  and $\eta_{j}$ is
the orthogonally scattered random measure in the spectral
representation of the stationary process $Y^s_j(\cdot)$.  Also by  the independence of
$\{U^j(\cdot)\}$ and $\{Y^s_j(\cdot)\}$, the cross spectral covariance
$F_{j,k}(d\lambda, d\lambda, t, u)=\big<\Phi_j(d\lambda,t), \Phi_k(d\lambda,u)\big>$ for $\lambda\neq\omega$  can be written as
$$F_{j,k}(d\lambda, d\lambda, t, u)
=\psi_{j-1}(t+u)G_{j-1}(d\lambda)I_{\{k,k+1\}}(j)+\psi_{j}(t+u)G_j(d\lambda)I_{\{k,k-1\}}(j),$$
where $\psi_{j}(t+u)=\big<U^j(t), U^j(u)\big>, $  $G_j(d\lambda)=E|\eta_{j}(d\lambda)|^2$, and  $F_{j,k}(d\lambda, d\omega, t, u)=0$.
\end{remark}

\begin{lemma}
Under the assumptions of  Remark 3.2-3.3, the cross covariance of the processes $X_j^p(t)$ and $Y^s_r(u)$ is bi-periodic,  with period $S$ in $t$ and $u$ and  period $T$ in $j$ and $r$, $\; \big<X_j^p(t), Y^s_r(u)\big>= \big<X_{j+T}^p(t), Y^s_{r+T}(u+S)\big>$,where $t \in B_j$ and $u\in B_r$ if and only if
\vspace{-1mm}
\be \big<X_j^p(t), Y^s_r(u)\big>=\frac{a_j^t}{a_j}\int_0^{2\pi}\int_{-\infty}^{\infty}
e^{i(j\lambda -u\omega)}\Delta_{r}(d\lambda, d\omega)\ee
where support of $\Delta_{r}(d\lambda, d\omega)\!=\! \big<\vartheta (d\lambda), \eta_{r}(d\omega)\big>$ is
$\Gamma \!=\{ (\lambda,\omega):\! \lambda T\!=S\omega-2\pi l, l \in {\mathbb Z } \}$.

\end{lemma}
{\bf Proof:}
By (3.9)  $\;<X_j^p(t), Y^s_r(u)>=\frac{a_j^t}{a_j}<X_j^p, Y^s_r(u)>$, by (3.11)
$\; X^p_j
=\int_{0}^{2\pi}e^{i\lambda j }\vartheta(d\lambda )$, and by (3.6)
$\;Y^s_r(u)
=\int_{-\infty}^{\infty}e^{i\omega u }\eta_{r}(d\omega)$, where $\eta_{r+T}(\cdot)=\eta_r(\cdot)$ . So we have (3.13), and one can easily verify that if $\, \Gamma \,$ is the support of the cross spectral covariance then  the covariance function has bi-periodic property.
On the other hand by assuming that the cross covariance is bi-periodic, using (3.9), equality of $\theta_r(j,u)=\theta_{r+T}(j+T,u+S)$ implies that for any natural number $N$, $$\; \theta_r(j,u)=\frac{1}{2N+1}\sum_{k=-N}^N \theta_{r+kT}(j+kT,u+kS)=\int_0^{2\pi}\!\!\!\int_{-\infty}^{\infty}\!\! e^{i(j\lambda-u\omega)}D_N(T\lambda-S\omega)\Delta_r(d\lambda,d\omega),$$ where $D_N(2\pi l)=1, \, l \in {\mathbb Z}$, and
$D_N(\sigma)=\frac{\sin (N+1/2)\sigma}{2(N+1/2)\sin(\sigma/2)}$ and converges to zero for $\sigma\neq 2\pi {\mathbb Z}$. This implies that the support of cross spectral density must be  contained in $\Gamma$.

\begin{theorem}
The spectral representation of the multi-component {\em PC-LS}
process $X(t)=X^{ls}(t)+X^p(t)$, where $X^p(t)$ and $X^{ls}(t)$ are
defined by $(\ref{vart-spec})$, $(\ref{sbn0-3})$ and  $(\ref{sbn0-3-1})$ can be written as
\be X(t)=\int_{-\infty}^{\infty}e^{i\lambda t}dZ_{m,t}(\lambda),\hspace{1cm}t\in B_m\label{sbn7},\ee
where  ${\bf dZ_{m,t}(\lambda)\!\! =\! \!\Phi_m(d\lambda,
t)+I_{[0,2\pi)}(\lambda)e^{i\lambda(m-t)}\bigg[\Upsilon(d\lambda,a_m^t)+A_{t,m}e^{-i\lambda}d \vartheta \bigg]}$, in which $\vartheta$, $\Phi_m$ and $\Upsilon$ are
defined by $(\ref{sbn5})$ and Remark 3.2.
The cross spectral covariance for $t\in B_m$, $u\in B_n$
{\bf\begin{eqnarray*}{\bf \big<dZ_{m,t}(\lambda), dZ_{n,u}(\omega)\big>}\!\!\!\!\!&=&\!\!\!{\bf F_{m,n}(d\lambda, d\omega, t, u)+K(m,n,\lambda,\omega)\bigg[ e^{-i(\lambda-\omega)}A_{t,m}A_{u,n}\tilde{\theta}(d\lambda,d\omega)}\\&&\hspace{-.5cm} {\bf +\Theta_{m,n}(d\lambda,d\omega,a_m^t,a_n^u)\bigg]\!\! +\! \frac{aa_n^u}{a_n}\Lambda_{m,n}(d\omega,d\lambda,t)\! +\! \frac{aa^t_m}{a_m}\Lambda_{n,m}(d\lambda,d\omega,u)}\\&&\hspace{.6cm}{\bf + A_{u,n}a \Lambda_{m,n-1}(d\omega,d\lambda,t) e^{-i\omega}+A_{t,m}a\Lambda_{n,m-1}(d\lambda,d\omega,u)e^{-i\lambda}}
\end{eqnarray*}}
where $A_{t,m}=1-a^t_m/a_m$,  $\,\Lambda_{l,k}(d\theta,d\eta,v)=e^{i(\theta k-v \eta)}\big(\Delta_l(d\theta,d\eta)+\Delta_{l-1}(d\theta,d\eta)\big)I_{[0,2\pi)}(d\theta ),\; l,k=m,n$, and   $\; \Theta,\;$  $\,
 F,\,$ and $\Delta$ are defined in {\em Remark 3.2},  {\em Remark 3.3}  and  {\em Lemma 3.3} respectively, and $K(m,n,\lambda,\omega)=e^{i\lambda(m-t)-i\omega(n-u)}I_{[0,2\pi)}(\lambda)I_{[0,2\pi)}(\omega)$.
\end{theorem}
{\bf Proof:}
By the result of Remarks 3.2-3.3  and relations (\ref{sbn0-3}),(\ref{sbn0-3-1}) we have
that $\{X(t), t\in {\mathbb R}^+\}$ has the time varying spectral representation  (\ref{sbn7}). Also by (3.9), Lemma 3.3, and Remark 3.2-3.3, the last assertion of the theorem can be easily obtained as the cross spectral covariance.

\end{document}